\documentclass{wstmp}
\usepackage{braket}
\usepackage{subfigure}

\newtheorem{program}{Program}[section]

\newcommand{\emptypartition}{\varnothing}

\newcommand{\NN}{\mathbb{N}}
\newcommand{\RR}{\mathbb{R}}

\newcommand{\AAA}{\mathcal{A}}
\newcommand{\BBB}{\mathcal{B}}
\newcommand{\PPP}{\mathcal{P}}

\newcommand{\FFF}{\mathcal{F}}
\newcommand{\VVV}{\mathcal{V}}

\newcommand{\codim}{\operatorname{codim}}

\newcommand{\partitionof}{\vdash}
\newcommand{\ideal}{\mathcal{I}}

\newcommand{\Stab}{\operatorname{Stab}}

\newcommand{\minelm}{\hat{0}}
\newcommand{\maxelm}{\hat{1}}
\newcommand{\numof}[1]{\left|#1\right|}

\newcommand{\figGraph}{%
\begin{picture}(100,40)
\put(50,30){\makebox(0,0){$\bullet$}}
\put(10,10){\makebox(0,0){$\bullet$}}
\put(20,10){\makebox(0,0){$\bullet$}}
\put(30,10){\makebox(0,0){$\bullet$}}
\put(50,10){\makebox(0,0){$\cdots$}}
\put(70,10){\makebox(0,0){$\bullet$}}
\put(80,10){\makebox(0,0){$\bullet$}}
\put(90,10){\makebox(0,0){$\bullet$}}
\put(50,35){\makebox(0,0)[b]{\tiny$1$}}
\put(10,5){\makebox(0,0)[t]{\tiny$2$}}
\put(20,5){\makebox(0,0)[t]{\tiny$3$}}
\put(30,5){\makebox(0,0)[t]{\tiny$4$}}
\put(90,5){\makebox(0,0)[t]{\tiny$l$}}
\thicklines
\put(35,10){\line(-1,0){25}}
\put(65,10){\line(1,0){25}}
\put(50,30){\line(-2,-1){40}}
\put(50,30){\line(-3,-2){30}}
\put(50,30){\line(-1,-1){20}}
\put(50,30){\line(1,-1){20}}
\put(50,30){\line(3,-2){30}}
\put(50,30){\line(2,-1){40}}
\end{picture}%
}

\newcommand{\scVi}{%
\begin{picture}(40,40)
\put(10,20){\makebox(0,0){$\bullet$}}
\put(30,20){\makebox(0,0){$\bullet$}}
\put(20,10){\makebox(0,0){$\bullet$}}
\put(20,30){\makebox(0,0){$\bullet$}}
\put(20,20){\makebox(0,0){$\bullet$}}
\put(15,5){\line(0,1){30}}
\put(15,35){\line(1,0){10}}
\put(15,5){\line(1,0){10}}
\put(25,5){\line(0,1){30}}
\put(5,15){\line(1,0){30}}
\put(35,15){\line(0,1){10}}
\put(5,15){\line(0,1){10}}
\put(5,25){\line(1,0){30}}
\thicklines
\put(30,20){\line(-1,1){10}}
\put(10,20){\line(1,-1){10}}
\put(10,20){\line(1,1){10}}
\end{picture}%
}
\newcommand{\scVix}{%
\begin{picture}(40,40)
\put(10,10){\makebox(0,0){$\bullet$}}
\put(10,30){\makebox(0,0){$\bullet$}}
\put(30,10){\makebox(0,0){$\bullet$}}
\put(30,30){\makebox(0,0){$\bullet$}}
\put(20,20){\makebox(0,0){$\bullet$}}
\put(5,5){\line(0,1){30}}
\put(5,35){\line(1,0){30}}
\put(5,5){\line(1,0){30}}
\put(35,5){\line(0,1){30}}
\end{picture}%
}
\newcommand{\scVviii}{%
\begin{picture}(40,50)
\put(10,10){\makebox(0,0){$\bullet$}}
\put(10,20){\makebox(0,0){$\bullet$}}
\put(10,30){\makebox(0,0){$\bullet$}}
\put(10,40){\makebox(0,0){$\bullet$}}
\put(30,30){\makebox(0,0){$\bullet$}}
\put(5,5){\line(0,1){40}}
\put(5,45){\line(1,0){10}}
\put(5,5){\line(1,0){10}}
\put(15,5){\line(0,1){40}}
\thicklines
\put(10,40){\line(2,-1){20}}
\put(10,10){\line(1,1){20}}
\put(10,20){\line(2,1){20}}
\put(10,30){\line(1,0){20}}
\end{picture}%
}
\newcommand{\scVvii}{%
\begin{picture}(40,40)
\put(10,10){\makebox(0,0){$\bullet$}}
\put(10,20){\makebox(0,0){$\bullet$}}
\put(10,30){\makebox(0,0){$\bullet$}}
\put(30,30){\makebox(0,0){$\bullet$}}
\put(30,10){\makebox(0,0){$\bullet$}}
\thicklines
\put(10,10){\line(1,0){20}}
\put(10,10){\line(1,1){20}}
\put(10,20){\line(2,1){20}}
\put(10,20){\line(2,-1){20}}
\put(10,30){\line(1,0){20}}
\put(10,30){\line(1,-1){20}}
\put(30,10){\line(0,1){20}}
\end{picture}%
}

\newcommand{\scVvi}{%
\begin{picture}(40,40)
\put(10,10){\makebox(0,0){$\bullet$}}
\put(10,20){\makebox(0,0){$\bullet$}}
\put(10,30){\makebox(0,0){$\bullet$}}
\put(30,30){\makebox(0,0){$\bullet$}}
\put(30,10){\makebox(0,0){$\bullet$}}
\put(5,5){\line(0,1){30}}
\put(5,35){\line(1,0){10}}
\put(5,5){\line(1,0){10}}
\put(15,5){\line(0,1){30}}
\thicklines
\put(10,10){\line(1,0){20}}
\put(10,10){\line(1,1){20}}
\put(10,20){\line(2,1){20}}
\put(10,30){\line(1,0){20}}
\put(30,10){\line(0,1){20}}
\end{picture}%
}
\newcommand{\scVv}{%
\begin{picture}(40,40)
\put(10,10){\makebox(0,0){$\bullet$}}
\put(10,20){\makebox(0,0){$\bullet$}}
\put(10,30){\makebox(0,0){$\bullet$}}
\put(30,10){\makebox(0,0){$\bullet$}}
\put(30,30){\makebox(0,0){$\bullet$}}
\thicklines
\put(10,10){\line(1,0){20}}
\put(10,20){\line(2,1){20}}
\put(10,20){\line(2,-1){20}}
\put(10,30){\line(1,0){20}}
\put(30,10){\line(0,1){20}}
\put(10,10){\line(0,1){20}}
\end{picture}%
}
\newcommand{\scViv}{%
\begin{picture}(40,40)
\put(10,10){\makebox(0,0){$\bullet$}}
\put(10,20){\makebox(0,0){$\bullet$}}
\put(10,30){\makebox(0,0){$\bullet$}}
\put(30,10){\makebox(0,0){$\bullet$}}
\put(30,30){\makebox(0,0){$\bullet$}}
\put(5,5){\line(0,1){30}}
\put(5,35){\line(1,0){10}}
\put(5,5){\line(1,0){10}}
\put(15,5){\line(0,1){30}}
\thicklines
\put(10,10){\line(1,0){20}}
\put(10,20){\line(2,1){20}}
\put(10,20){\line(2,-1){20}}
\put(10,30){\line(1,0){20}}
\put(30,10){\line(0,1){20}}
\end{picture}%
}
\newcommand{\scViii}{%
\begin{picture}(40,40)
\put(10,20){\makebox(0,0){$\bullet$}}
\put(30,20){\makebox(0,0){$\bullet$}}
\put(20,10){\makebox(0,0){$\bullet$}}
\put(20,30){\makebox(0,0){$\bullet$}}
\put(20,20){\makebox(0,0){$\bullet$}}
\put(15,5){\line(0,1){30}}
\put(15,35){\line(1,0){10}}
\put(15,5){\line(1,0){10}}
\put(25,5){\line(0,1){30}}
\thicklines
\put(30,20){\line(-1,0){10}}
\put(30,20){\line(-1,-1){10}}
\put(30,20){\line(-1,1){10}}
\put(10,20){\line(1,-1){10}}
\put(10,20){\line(1,1){10}}
\end{picture}%
}

\newcommand{\scVii}{%
\begin{picture}(40,40)
\put(10,20){\makebox(0,0){$\bullet$}}
\put(30,20){\makebox(0,0){$\bullet$}}
\put(20,10){\makebox(0,0){$\bullet$}}
\put(20,30){\makebox(0,0){$\bullet$}}
\put(20,20){\makebox(0,0){$\bullet$}}
\thicklines
\put(20,20){\line(1,0){10}}
\put(30,20){\line(-1,-1){10}}
\put(30,20){\line(-1,1){10}}
\put(20,10){\line(0,1){20}}
\put(10,20){\line(1,-1){10}}
\put(10,20){\line(1,1){10}}
\end{picture}%
}

\newcommand{\scUniv}{%
\begin{picture}(180,40)
\put(10,10){\makebox(0,0){$\bullet$}}
\put(30,30){\makebox(0,0){$\bullet$}}
\put(50,10){\makebox(0,0){$\bullet$}}
\put(70,30){\makebox(0,0){$\bullet$}}
\put(90,10){\makebox(0,0){$\bullet$}}
\put(110,30){\makebox(0,0){$\bullet$}}
\put(130,10){\makebox(0,0){$\bullet$}}
\put(155,20){\makebox(0,0){$\cdots$}}
\put(15,10){\makebox(0,0)[l]{$1$}}
\put(35,30){\makebox(0,0)[l]{$2$}}
\put(55,10){\makebox(0,0)[l]{$3$}}
\put(75,30){\makebox(0,0)[l]{$4$}}
\put(95,10){\makebox(0,0)[l]{$5$}}
\put(115,30){\makebox(0,0)[l]{$6$}}
\put(135,10){\makebox(0,0)[l]{$7$}}
\put(5,5){\line(1,0){160}}
\put(5,5){\line(0,1){10}}
\put(165,5){\line(0,1){10}}
\put(5,15){\line(1,0){160}}
\put(25,25){\line(1,0){150}}
\put(25,25){\line(0,1){10}}
\put(175,25){\line(0,1){10}}
\put(25,35){\line(1,0){150}}
\thicklines
\put(10,10){\line(1,1){20}}
\put(50,10){\line(1,1){20}}
\put(50,10){\line(-1,1){20}}
\put(90,10){\line(1,1){20}}
\put(90,10){\line(-1,1){20}}
\put(130,10){\line(1,1){12}}
\put(130,10){\line(-1,1){20}}
\end{picture}%
}

\newcommand{\scIVi}{%
\begin{picture}(40,40)
\put(10,20){\makebox(0,0){$\bullet$}}
\put(30,20){\makebox(0,0){$\bullet$}}
\put(20,10){\makebox(0,0){$\bullet$}}
\put(20,30){\makebox(0,0){$\bullet$}}
\thicklines
\put(30,20){\line(-1,0){20}}
\put(30,20){\line(-1,-1){10}}
 \put(30,20){\line(-1,1){10}}
\put(10,20){\line(1,0){20}}
\put(10,20){\line(1,-1){10}}
\put(10,20){\line(1,1){10}}
\end{picture}%
}

\newcommand{\scIViii}{%
\begin{picture}(40,40)
\put(10,20){\makebox(0,0){$\bullet$}}
\put(30,20){\makebox(0,0){$\bullet$}}
\put(20,10){\makebox(0,0){$\bullet$}}
\put(20,30){\makebox(0,0){$\bullet$}}
\put(20,35){\line(-1,-1){15}}
\put(20,5){\line(-1,1){15}}
\put(20,5){\line(1,0){3}}
\put(20,35){\line(1,0){3}}
\put(23,5){\line(0,1){30}}
\thicklines
\put(30,20){\line(-1,0){20}}
\put(30,20){\line(-1,-1){10}}
 \put(30,20){\line(-1,1){10}}
\end{picture}%
}
\newcommand{\scIVii}{%
\begin{picture}(40,40)
\put(10,20){\makebox(0,0){$\bullet$}}
\put(30,20){\makebox(0,0){$\bullet$}}
\put(20,10){\makebox(0,0){$\bullet$}}
\put(20,30){\makebox(0,0){$\bullet$}}
\put(20,35){\line(-1,-1){15}}
\put(20,35){\line(1,-1){15}}
\put(20,5){\line(-1,1){15}}
\put(20,5){\line(1,1){15}}
\end{picture}%
}

\newcommand{\scIIIi}{%
\begin{picture}(40,40)
\put(10,10){\makebox(0,0){$\bullet$}}
\put(30,10){\makebox(0,0){$\bullet$}}
\put(20,30){\makebox(0,0){$\bullet$}}
\thicklines
\put(10,10){\line(1,2){10}}
\put(30,10){\line(-1,2){10}}
\put(10,10){\line(1,0){20}}
\end{picture}%
}

\newcommand{\scIIIii}{%
\begin{picture}(40,40)
\put(10,10){\makebox(0,0){$\bullet$}}
\put(30,10){\makebox(0,0){$\bullet$}}
\put(20,30){\makebox(0,0){$\bullet$}}
\put(5,7){\line(1,2){15}}
\put(35,7){\line(-1,2){15}}
\put(5,7){\line(1,0){30}}
\end{picture}%
}

\newcommand{\runningexampleI}{%
\begin{picture}(70,70)
\put(26.66,26.66){\makebox(0,0){$\bullet$}}
\put(10,10){\makebox(0,0){$\bullet$}}
\put(10,35){\makebox(0,0){$\bullet$}}
\put(35,10){\makebox(0,0){$\bullet$}}
\end{picture}%
}
\newcommand{\runningexampleII}{%
\begin{picture}(70,70)
\put(0,10){\line(1,0){70}}
\put(10,0){\line(0,1){70}}
\put(0,0){\line(1,1){50}}
\put(0,40){\line(2,-1){70}}
\put(40,0){\line(-1,2){35}}
\put(45,0){\line(-1,1){45}}
\end{picture}%
}
\newcommand{\runningexampleIII}{%
\begin{picture}(70,70)
\put(0,10){\line(1,0){70}}
\put(10,0){\line(0,1){70}}
\put(0,0){\line(1,1){50}}
\put(0,40){\line(2,-1){70}}
\put(40,0){\line(-1,2){35}}
\put(45,0){\line(-1,1){45}}
\put(26.66,26.66){\makebox(0,0){$\bullet$}}
\put(10,10){\makebox(0,0){$\bullet$}}
\put(10,35){\makebox(0,0){$\bullet$}}
\put(35,10){\makebox(0,0){$\bullet$}}
\put(10,60){\makebox(0,0){$\circ$}}
\put(60,10){\makebox(0,0){$\circ$}}
\put(22.5,22.5){\makebox(0,0){$\circ$}}
\end{picture}%
}

\begin{document}

\title{On computation of the characteristic polynomials of the discriminantal arrangements and the arrangements generated by generic points} 

\author{Yasuhide NUMATA}
\address{Mathematical Informatics, 
Graduate School of Information Science and Technology, 
University of Tokyo.\\
JST CREST.}

\author{Akimichi TAKEMURA}
\address{Mathematical Informatics, 
Graduate School of Information Science and Technology, 
University of Tokyo.\\
JST CREST.}

\begin{abstract}
In this article we give a computational study of combinatorics of the discriminantal arrangements. The discriminantal arrangements are parametrized by two positive integers $n$ and $k$ such that $n>k$.
The intersection lattice of the discriminantal arrangement with the parameter $(n,k)$
is isomorphic to the intersection lattice of
the hyperplane arrangement generated by $n$ generic points in the $d$-dimensional
vector space where $d=n-k-1$. The combinatorics of the discriminantal arrangements is very hard,
except for the special cases of 
the Boolean arrangements ($k=0$) and  the braid arrangements ($k=1$).
We review some results on the intersection lattices
of the arrangements generated by generic points 
and use them to obtain some computational results
on the characteristic polynomials of the discriminantal arrangements.
\end{abstract}


\bodymatter

\section{Introduction}
\label{sec:intro}
The main subject of this article is the combinatorics 
of the discriminantal arrangement
introduced by Manin and Schechtman \cite{manin-schechtman1989}. 
The discriminantal arrangement 
$\BBB(n,k)$, $n>k$, 
is defined as follows (cf.\ Section 5.6 of Orlik and Terao\cite{orlik-terao-book}):
Let $K$ be a field of characteristic zero
and 
let 
$\AAA_0=\Set{H_1,\ldots,H_{n}}$
a general position arrangement of $n$ hyperplanes
in a $k$-dimensional vector space over $K$.
Fix a normal vector $\phi_i$ for each hyperplane $H_i$.
For $a \in K^n$, 
we define $\AAA_{a}$  to be the hyperplane arrangement 
$\AAA_{a}=\Set{H_1+a_1 \phi_1,\ldots,H_{n}+a_{n}\phi_{n}}$
obtained from $\AAA_0$ by parallel translation.
The set
of the parameter $a\in K^n$ 
such that the hyperplanes of $\AAA_{a}$ are in a general position
is a complement  
$K^{n}\setminus \BBB$ of some central hyperplane
arrangement $\BBB$ in $K^{n}$. 

When $\AAA_0$ is generic in the sense of 
of Athanasiadis\cite{athanasiadis2000},   
the combinatorial structure of $\BBB$
depends only on $n$ and $k$, and
is called the discriminantal arrangement $\BBB(n,k)$. 
Let 
\[
d=n-k-1.
\]
Bayer and Brandt\cite{bayer-brandt1997} conjectured that 
the intersection lattice of 
$\BBB(n,k)$ is isomorphic 
to a lattice 
$L_{n,d}$ defined in set-theoretical terminology
(see Section \ref{sec:def} for the definition of $L_{n,d}$) and 
Athanasiadis\cite{athanasiadis2000} proved the conjecture.

Falk\cite{falk1994} showed the equivalence between  the discriminantal arrangements
and the  arrangements generated by generic points.
Regard $\phi_i^t$ as $k$-dimensional row vectors
and consider the image $P$ of the linear map 
from $K^{k}$ to 
$K^{n}$ defined by the $n\times k$-matrix
\begin{align*}
\begin{pmatrix}
\phi_1^t\\
\vdots\\
\phi_{n}^t
\end{pmatrix}.
\end{align*}
In other words, $P$ is the subspace of $K^n$ 
spanned by the column vectors of the matrix.
The image of the hyperplane $H_i$ is $\Set{x\in P|x_i=0}$.
Since we assume that $\phi_1,\ldots, 
\phi_n$  are in general position,
the dimension of the subspace $P$ is $k$.
Consider the orthogonal complement of $P$ 
and fix a basis $\Set{v_1,\ldots,v_{d+1}}$ of the complement.
Let $w_i^t$ be the $i$-th $(d+1)$-dimensional row vector of
the $n\times (d+1)$-matrix $(v_1,\ldots,v_{d+1})$.
Falk showed that 
the $n$ points $w_1,\ldots, w_n\in K^{d+1}$ 
are in general position, and that
the essential part of the discriminantal arrangement $\BBB(n,k)$ 
is equivalent to the central hyperplane arrangement
consisting of 
all hyperplanes in $K^{d+1}$
spanned by subsets of $\Set{w_1,\ldots, w_{n}}$ 
of size $d$. 
It is easy to see that
the latter arrangement
can be obtained 
by coning without the additional coordinate hyperplane  
from 
the non-central hyperplane arrangement $\AAA_{n,d}$ 
generated by generic 
$n$ points in $d$-dimensional vector space.
(See Section \ref{sec:def} for the definition of $\AAA_{n,d}$.) 



Recently 
the present authors showed that
the order ideals  of the intersection lattice of $\AAA_{n,d}$ 
can be decomposed into
direct products of smaller lattices corresponding to smaller dimensions\cite{1009.3676}.
The decomposition implies some identities of the M\"obius functions 
and the characteristic polynomials of the intersection lattices.
Moreover they give a way to compute 
the M\"obius functions 
and the characteristic polynomials of the intersection lattices
as polynomials in $n$. 
We first review this result and  then discuss computation of
the characteristic polynomials of the discriminantal arrangements.

The organization of this article is the following.
In Section \ref{sec:def}, we set up our definitions and notation.  
In Section \ref{sec:computation}, we recall the definition of 
the M\"obius function 
and the characteristic polynomial
of an intersection lattice.
These two sections are largely based on our previous paper\cite{1009.3676}.
Then we give some examples of computational results in Section \ref{sec:tables}.
Fixing $k$, we have a family of hyperplane arrangement  $\AAA_{d+k+1,d}$
parametrized by $d$
(or equivalently $\BBB(n,k)$, $n=d+k+1$, parametrized by $n$). 
In Section \ref{sec:disc}, 
we discuss the structure of the intersection lattices
of hyperplane arrangements of the family for $k=0,1,2$.
In Section \ref{sec:computation:disc}, 
we present some computational results of the characteristic polynomials
for $\AAA_{d+3,d}$.

\section{Intersection lattice and order ideal}
\label{sec:def}
In this section, we define the hyperplane arrangement $\AAA_{n,d}$ 
generated by generic points
and recall an interpretation
of the intersection lattice of $\AAA_{n,d}$
in set-theoretical terminology.
We also state some fundamental structure of 
order ideals of the intersection lattice of $\AAA_{n,d}$. 

In our previous paper\cite{1009.3676} we mainly considered fixed $d$ and arbitrary $n>d$.
For the case of fixed $k$ and arbitrary $d$, it is more convenient to write
$n$ as $d+k+1$. In the following we either write $\AAA_{n,d}$ or
$\AAA_{d+k+1,d}$ depending on the context.

Let 
$\PPP=\Set{p_1,\ldots,p_{n}}$ be
a set of $n$ points 
in $V=K^d$. 
We assume 
that $p_1,\ldots,p_{n}$ are generic 
in the sense of Athanasiadis\cite{athanasiadis2000}.  
For $X\subset \PPP$, define $H_{X}$ to be 
the affine hull of $X$.
Let 
$\AAA=\Set{H_{X}| X\subset \PPP, \numof{X} = d}$
be the set of all affine hyperplanes defined by $H_{X}$ for some $X\subset\PPP$, 
$\numof{X}=d$.
Since we consider generic points,
the combinatorial properties of the arrangement $\AAA$ 
do not depend on the points. 
We define $\AAA_{n,d}$ to be the arrangement $\AAA$,
and $L(\AAA_{n,d})$ to be the intersection lattice of
$\AAA_{n,d}$, 
i.e., the set
$\set{H_1 \cap \cdots \cap H_l | H_1,\ldots, H_l \in \AAA_{n,d}}$,
ordered by reverse inclusion.  
Contrary to the usual convention, 
here we consider that 
$\emptyset = \bigcap_{X\colon \numof{X}=d} H_X$ belongs to $L(\AAA_{n,d})$, 
so that 
$L(\AAA_{n,d})$ is not only a poset but also a lattice\cite{stanley-introduction}.
In the usual convention, 
this corresponds to the coning $c\AAA_{n,d}$ of $\AAA_{n,d}$, 
except that we do not add a coordinate hyperplane.  
The reason for this unconventional definition
is that $\emptyset \in L(\AAA_{n,d})$ plays an essential role
for recursive description of $L(\AAA_{n,d})$.

\begin{example}
Consider the two-dimensional vector space $\RR^2$.
Let $\PPP$ be the set of points in Figure 1(a). 
In this case,
$\AAA$ is the set of the lines in  Figure 1(b). 
The set of points, i.e., the elements of codimension two, in the intersection lattice
consists of seven points in  Figure 1(c). 
The four black points in  Figure 1(c) 
are original points in $\PPP$.
On the other hand,
the three white points in  Figure 1(c) 
are 
new points described as the intersection of two lines.
\begin{figure}[t]
\hfill
\subfigure[$\PPP$]{\runningexampleI}
\label{fig:1-1}
\hfill
\subfigure[$\AAA$]{\runningexampleII}
\label{fig:1-2}
\hfill
\subfigure[$L(\AAA_{4,2})$]{\runningexampleIII}
\label{fig:1-3}
\hfill{}
\caption{Example of $L(\AAA_{4,2})$}
\end{figure}
\end{example}

For the rest of this section we write $n$ as $d+k+1$.
We  recall an interpretation of 
$L(\AAA_{d+k+1,d})$ in set-theoretical terminology\cite{athanasiadis2000}.
\begin{definition}
For a finite set $X$, we define
\begin{align*}
\codim^k (X) &=\numof{X}-k.
\end{align*}
For distinct finite sets $S_1,\ldots,S_l$, we define
\begin{align*}
\rho^k(\Set{S_1,\ldots,S_l})
&=\codim^k (S_1)+\cdots+\codim^k (S_l), \\
D^k(\Set{S_1,\ldots,S_l})
&=
\codim^k (\bigcup_i S_i) -\rho^k(\Set{S_1,\ldots,S_l}).
\end{align*}
We define
$\rho^k(\emptyset)=\rho^k(\Set{})=0$.  
\end{definition}
\begin{definition}
\label{def:codim-new}
For $d>0$ and $k\geq 0$, we define $L^{d+k+1,d}$ by
\begin{align*}
L^{d+k+1,d}=
\Set{S \subset 2^{[d+k+1]} | 
\begin{array}{c}
\text{$D^k(S')>0$ for all $S'\subset S$ with $\numof{S'}>1$.}\\
\text{$k+1 \leq \numof{S_i}\leq d+k+1$ for all $S_i \in S$}.
\end{array}
}
\end{align*}
with ordering defined by
\begin{align*}
S<S' \iff 
\begin{cases}
\rho^k(S)<\rho^k(S')
& \text{and}\\
\text{$\forall S_i \in S$,
$\exists S_j' \in S'$ such that 
$S_i \subset S'_j$.}
\end{cases}
\end{align*}
\end{definition}

\begin{remark}
The poset $L^{d+k+1,d}$ 
is ranked by the rank function $\rho^k$.
The minimum element of $L^{d+k+1,d}$
is $\emptyset$,
and 
the maximum element of $L^{d+k+1,d}$ is
$\Set{\set{1,\ldots,d+k+1}}$
of rank $d+1$.
We use the symbol $\minelm$ (resp.\ $\maxelm$) to denote 
the minimum (resp.\ maximum) element 
of $L^{d+k+1,d}$.
If $d=0$, then
the poset
$L^{k+1,0}$ consists of two elements $\set{\minelm, \maxelm}$
for any $k\geq 0$.
\end{remark}

Fix $d$ and $k$. 
Then we have a bijection from $2^{[d+k+1]}$ to $2^{[d+k+1]}$ 
which maps 
a set $S_l$ to its complement $T_l=[d+k+1]\setminus S_l$.
Hence the definition of 
$\codim^k$,
$\rho^k$,
$D^k$
and
$L^{d+k+1,d}$
can be rewritten with the complements as follows:
\begin{definition}
\label{def:codim-old}
For a finite set  $X$, 
we define 
\begin{align*}
\codim_d(X) = d+1-\numof{X}.
\end{align*}
For distinct finite sets $T_1,\ldots, T_l$,
we define 
\begin{align*}
\rho_d(\Set{T_1,\ldots,T_l})
&= \codim_d (T_1) +\cdots+\codim_d (T_l),\\
D_d(\Set{T_1,\ldots,T_l})
&=\codim_d(T_1\cap\cdots\cap T_l)-\rho_d(\Set{T_1,\ldots,T_l}).
\end{align*}
We define $\rho_d(\emptyset)=\rho_d(\Set{ }) = 0$.
\end{definition}
This is the definition we employed in our previous paper\cite{1009.3676}.
In this article we need to use both Definition \ref{def:codim-new}
and Definition \ref{def:codim-old}. We distinguish them by 
superscripts and subscripts.

\begin{remark}
Let $T_1,\ldots,T_l \subset [d+k+1]$ and $S_i=[d+k+1] \setminus T_i$.
Then 
\begin{align*}
\codim_d (T_i) 
&=d+1 -\numof{T_i}\\ 
&=d+1+\numof{S_i}-(d+k+1)\\
&=\numof{S_i}-k\\
&=\codim^k(S_i).
\end{align*}
We also have
\begin{align*}
\codim_d(\bigcap_i T_i)=\codim^k(\bigcup_i S_i),
\end{align*}
since $\bigcap_i T_i=[d+k+1]\setminus \bigcup_i S_i$.
Hence we have
$\rho_d(T)=\rho^k(S)$ and
$D_d(T)=D^k(S)$ for $T=\Set{T_1,\ldots,T_l}$
and $S=\Set{S_1,\ldots,S_l}$.
\end{remark}

\begin{definition}
\label{def:DefOfLattice}
For $d>0$ and $k\geq 0$, 
we define $L_{d+k+1,d}$ by
\begin{align*}
L_{d+k+1,d}=
\Set{T\subset 2^{\Set{1,\ldots,d+k+1 }} | 
\begin{array}{c}
\text{$D_d(T')>0$ for all $T'\subset T$ with $\numof{T'}>1$.}\\
\text{$0\leq \numof{T_i} \leq d$ for all $T_i \in T$.}
\end{array}
}
\end{align*}
with ordering  $<$ defined by
\begin{align*}
T<T' 
&\iff 
\begin{cases}
\rho_d(T) < \rho_d(T') & \text{and} \\
\text{$\forall T_i \in T$, $\exists T_j' \in T' $ such that $T_j' \subset T_i$.}
\end{cases}
\end{align*} 
\end{definition}
The map 
$L^{d+k+1,d} \ni S \mapsto \Set{[d+k+1] \setminus S_i | S_i \in S} \in L_{d+k+1,d}$ 
is an order-preserving bijection.
Hence $L^{d+k+1,d}$ is isomorphic to  $L_{d+k+1,d}$
as posets.


We now state some fundamental structures of the intersection lattice $L_{d+k+1,d}$.
For $T\in L_{d+k+1,d}$,
we define $\ideal_{d+k+1,d}(T)$ to be the {\it order ideal} generated by $T$,
that is, $\ideal_{d+k+1,d}(T)=\Set{T' \in L_{d+k+1,d}|T'\leq T}$.
Then we have the following theorems\cite{1009.3676}.
\begin{theorem}
\label{lem:single-set-ideal}
For $\set{T_1}\in L_{d+k+1,d}$, 
$\ideal_{d+k+1,d}(\{T_1\})$ is isomorphic to 
$L_{d+k+1-\numof{T_1},d-\numof{T_1}}=L_{k+\codim_d(T_1),\codim_d(T_1)-1}$
as posets.
\end{theorem}
\begin{theorem}
\label{thm:structure of ideal is directprod}
Let $T\in L_{d+k+1,d}$, $T\neq \minelm$. 
Then the ideal $\ideal_{d+k+1,d}(T)$
is isomorphic to the direct product 
$\prod_{T_i \in T} \ideal_{d+k+1,d}(\{T_i\})$
as posets.  They are also isomorphic to 
$\prod_{T_i \in T} L_{k+\codim_d(T_i),\codim_d(T_i)-1}$.
\end{theorem}

\section{Computation of the M\"obius function and the characteristic polynomial}
\label{sec:computation}

In this section 
we apply Theorem \ref{thm:structure of ideal is directprod} 
to the M\"obius function and the characteristic polynomial 
of the intersection lattice $L_{n,d}$. 

\subsection{M\"obius function and characteristic polynomial}
\label{subsec:mobius}
By Theorem \ref{thm:structure of ideal is directprod},
the combinatorial structure of $\ideal_{n,d}(T)$ depends 
only on $\codim(T_i)$.
Therefore we introduce the notion of the type of $T$.

Let $d$ be a nonnegative integer.
We call 
a weakly-decreasing sequence $\delta=(\delta_1,\delta_2,\ldots)$
of nonnegative integers such that $\sum_{i} \delta_i=d$
a {\em partition} of $d$.
We allow one or more zeros to occur at the end,
or equivalently,
$(\delta_1,\delta_2,\ldots,\delta_l)=(\delta_1,\delta_2,\ldots,\delta_l,0)$
for any partition $\delta=(\delta_1,\delta_2,\ldots,\delta_l)$.
We write $\delta\partitionof d$ to say that
$\delta$ is a partition of $d$.
For example, $\Set{ \delta \partitionof 3 }=\set{(3),(2,1),(1,1,1)}$,
and $\Set{ \delta \partitionof 0 }$ 
is the set consisting of  the unique partition of zero, 
which is denoted by $\emptypartition$.

\begin{definition}
For $T=\set{T_1, \dots, T_l}\in L_{n,d}$
with $\numof{T_1}\leq \dots \leq \numof{T_l}$,
we call 
$\gamma_d(T)=(\codim_d(T_1),\dots,\codim_d(T_l)) \partitionof \rho_d(T)$
the {\em type} of $T$.  
We define $\Gamma(n,d)=\Set{\gamma_d(T)|T\in L_{n,d}}$.
\end{definition}
\begin{example}
For each $d$, 
$\gamma_d(\minelm)=\emptypartition$ and $\gamma_d(\maxelm)=(d+1)$.
\end{example}
\begin{remark}
If $d$ is fixed and $n$ is sufficiently large,  
then 
\begin{align*}
\Gamma(n,d)=\Set{\gamma\partitionof i|i=0,1,\ldots,d}\cup \Set{(d+1)}.
\end{align*}
\end{remark}

The {\it M\"obius function} $\mu_{n,d}$ 
and the {\it characteristic polynomial} $\chi_{n,d}(t)$
of the poset $L_{n,d}$ is defined as usual\cite{stanley-introduction}:
\begin{align*}
 &\mu_{n,d}(T,T)=1,   \quad 
\sum_{T''\colon T\leq T''\leq T'} \mu_{n,d}(T,T'')=0 \quad (T < T'),\\
&\mu_{n,d}(T) = \mu_{n,d}(\minelm,T),\\
&\chi_{n,d}(t) =\sum_{T\in L_{n,d}} \mu_{n,d}(T) t^{d+1-\rho_d(T)}.
\end{align*}

The usual characteristic polynomial $\chi(\AAA_{n,d},t)$
of the non-central arrangement $\AAA_{n,d}$ 
is given as
\begin{align*}
\chi(\AAA_{n,d},t)
&=\sum_{T\in L_{n,d}, \, T\neq \maxelm} \mu_{n,d}(T) t^{d-\rho_d(T)}
=\frac{\chi_{n,d}(t)- \mu_{n,d}(\maxelm)}{t}.
\end{align*}
Since  $\chi_{n,d}(1)=0$, 
\begin{align}
\label{eq:mumaxphi}
  \mu_{n,d}(\maxelm)
=- \sum_{T\in L_{n,d}, \, T\neq \maxelm} \mu_{n,d}(T).
\end{align}

Next we evaluate the value of M\"obius function of $L_{n,d}$.
The M\"obius function for the direct product of posets
is written as the product of the M\"obius functions of posets.
Hence 
Theorem \ref{thm:structure of ideal is directprod} implies
the following theorem.

\begin{theorem}
\label{thm:decomposition}
For $T\in L_{n,d}$ and $T\neq \minelm$, we have
\begin{align*}
\mu_{n,d}(T)=
\prod_{T_i \in T}
\mu_{d+k+1-\numof{T_i},d-\numof{T_i}}(\maxelm)=
\prod_{T_i \in T}
\mu_{k+\codim_{d}(T_i),\codim_{d}(T_i)-1}(\maxelm).
\end{align*}
\end{theorem}
\begin{remark} 
Since $L_{k+1,0}=\set{\minelm,\maxelm}$,
we have  $\mu_{k+1,0}(\minelm)=1$ and $\mu_{k+1,0}(\maxelm)=-1$.
The equation $\mu_{k+1,0}(\maxelm)=-1$
implies that
$\mu_{k+\codim_{d}(T_i),\codim_{d}(T_i)-1}(\maxelm)=-1$
if $\codim_{d}(T_i)=1$.
\end{remark}

The value of the M\"obius function 
depends only on 
the type  $\gamma_d(T)$ of $T$.
Therefore, from now on 
we denote $\mu_{n,d}(T)=\mu_{n,d}(\gamma_d(T))$.
Then we can rewrite Theorem \ref{thm:decomposition} as
follows.
\begin{corollary}
\label{cor:mobius:decomp}
For $\gamma=(\gamma_1,\ldots,\gamma_l)\in \Gamma(n,d) \setminus\Set{\emptypartition}$,
\begin{align*}
\mu_{n,d}(\gamma)=
\prod_{i=1}^{l}\mu_{n-d+\gamma_i-1,\gamma_i-1}(\maxelm).
\end{align*}
\end{corollary}

\begin{corollary}
\label{cor:type}
If $d'\leq d$, then
\begin{align*}
\mu_{d+k+1,d}(\gamma)
=
\mu_{d'+k+1,d'}(\gamma)
\end{align*}
for $\gamma\in \Gamma(d'+k+1,d')$.
\end{corollary}

\subsection{Number of elements of the intersection lattice}
\label{subsec:number-of-elements}

Define $L_{n,d}(\gamma)$ and $\lambda_{n,d}(\gamma)$
by
\[
L_{n,d}(\gamma)=\Set{T \in L_{n,d}| \gamma_d(T)=\gamma},\quad 
\lambda_{n,d}(\gamma) = \numof{L_{n,d}(\gamma)}
\]
for $\gamma\in \Gamma(n,d)$.
Then from the results in the previous sections we have
\begin{align*}
\chi_{n,d}(t) 
&=
 \mu_{n,d}(\maxelm)+
\sum_{i=0}^{d}
\sum_{\gamma \partitionof i} 
\lambda_{n,d}(\gamma)\mu_{n,d}(\gamma) t^{d+1-i}
, 
\\
  \mu_{n,d}(\maxelm)
&=- 
\sum_{i=0}^{d}
\sum_{\gamma \partitionof i} 
 \lambda_{n,d}(\gamma)\mu_{n,d}(\gamma).
\end{align*}
Here we discuss how to compute the number $\lambda_{n,d}(\gamma)$.

From the viewpoint of computation,
tuples are easier than sets. 
Therefore we define $\tilde L^l_{n,d}$ and $\tilde L_{n,d}(\gamma)$ by
\begin{align*}
\tilde L^l_{n,d}&=
\Set{(T_1,\ldots,T_l)\in
 \left(2^{\Set{1,\ldots,n}}\right)^l
|
\begin{array}{c}
\Set{T_1,\ldots,T_l} \in L_{n,d}.\\
T_i\neq T_j (i\neq j).
\end{array}
},\\
\tilde L_{n,d}(\gamma)&=
\Set{
(T_1,\ldots,T_l)\in \tilde L^l_{n,d}
|
\text{$\codim_d(T_i) = \gamma_i$ for each $i$.}
}
\end{align*}
for $\gamma=(\gamma_1,\ldots,\gamma_l)\in \Gamma(n,d)$.
Obviously we have
 \begin{align*}
 \numof{\tilde L_{n,d}(\gamma)}
& =\lambda_{n,d}(\gamma)\cdot\numof{\Stab_{\mathfrak{S}_l}(\gamma)},
 \end{align*}
where
$\Stab_{\mathfrak{S}_l}(\gamma)$ stands for
the stabilizer of the symmetric group $\mathfrak{S}_l$ fixing 
a partition 
$\gamma=(\gamma_1, \dots, \gamma_l)$,
i.e.,
$\Set{\sigma\in \mathfrak{S}_l|\text{$\gamma_i=\gamma_{\sigma(i)}$ for each $i$}}$.

Each element $T=(T_1,\ldots,T_l)\in \tilde L^l_{n,d}$ is characterized by
information on `symbols' in each $T_i$.
Conversely, 
information on indices $i$ such that each `symbol' belongs to $T_i$
can  characterize $T$.
For 
 $T=(T_1,\ldots,T_l)\in\tilde L_{n,d}(\gamma)$ 
 and a subset $I \subset \Set{1,\ldots,l}$,
 define $\tau_T(I)$ by
\begin{align*}
 \tau_T(I)&=\Set{t| t\in T_i \iff i \in I}
=\bigcap_{i\in I} T_i \setminus \bigcup_{i\not\in I} T_i .
\end{align*}
For 
each $T=(T_1,\ldots,T_l)\in\tilde L^l_{n,d}$,
$\tau_T$ is
a map from $2^{\set{1,\ldots,l}}$ to $2^{\Set{1,\ldots,n}}$.
The map 
 $\tilde L^l_{n,d}\ni T \mapsto \tau_T \in \Set{ \pi \colon  2^{\set{1,\ldots,l}} \to 2^{\Set{1,\ldots,n}}}$ 
is an injection.
Let $\tilde N(n,d;\gamma)=\Set{\tau_T |T \in \tilde  L_{n,d}(\gamma)}$.
Then $\numof{\tilde N(n,d;\gamma)}=\numof{\tilde L^l_{n,d}(\gamma)}$.

For $\tau \in N(n,d;\gamma)$,
define a map
\begin{align*}
\bar \tau_T \colon 
2^{\set{1,\ldots,l}} \ni I \mapsto \numof{\tau_T(I)} \in \NN
\end{align*}
and let 
\begin{align*}
N(n,d;\gamma)=
\Set{
\bar \tau_T \colon 2^{\set{1,\ldots,l}} \to \NN
|
\tau_T \in \tilde N(n,d;\gamma)
}.
\end{align*}
Moreover, for $\nu \in N(n,d;\gamma)$, define 
$\tilde N(n,d;\gamma,\nu)$ by
\begin{align*}
\tilde N(n,d;\gamma,\nu)=
\Set{\tau_T\in \tilde N(n,d;\gamma)|
\text{$\numof{\tau_T(I)}= \nu(I)$ for each $I \in 2^{\Set{1,\ldots,l}}$.} 
}.
\end{align*}
Then 
$\tilde N(n,d;\gamma)=\coprod_{\nu \in N(n,d;\gamma)} \tilde N(n,d;\gamma,\nu)$,
which implies
\begin{align*}
\lambda_{n,d}(\gamma)\cdot\numof{\Stab_{\mathfrak{S}_l}(\gamma)}
=\numof{\coprod_{\nu \in N(n,d;\gamma)} \tilde N(n,d;\gamma,\nu)}.
\end{align*}
It is not difficult to see that
\begin{align*}
 \numof{\tilde N(n,d;\gamma,\nu)}
& =\binom{n}{\nu(I_1), \nu(I_2), \ldots ,\nu(I_{2^l})},
\end{align*}
where 
$2^{\Set{1,\ldots,l}}=\Set{I_1,\ldots,I_{2^l}}$
and
$\binom{n_1+\cdots+n_m}{n_1,n_2,\ldots,n_m}$ stands for 
the multinomial coefficient $(n_1+\cdots+n_m)!/(n_1!\,n_2!\,\cdots n_m!)$.
Therefore we have the following proposition.

\begin{proposition}
\label{prop:number-type}
For $\gamma\in \Gamma(n,d) \setminus\Set{\emptypartition}$,
$\lambda_{n,d}(\gamma)$ is given by
\begin{align}
\lambda_{n,d}(\gamma)
&=\frac{1}{\prod_{s=1}^d m_s(\gamma) !}
\sum_{\nu\in N(n,d;\gamma)} 
\binom{n}{\nu(I_1), \nu(I_2),\ldots, \nu(I_{2^l})},
\label{eq:number-type}
\end{align}
where 
$2^{\Set{1,\ldots,l}}=\Set{I_1,\ldots,I_{2^l}}$
and
$m_s(\gamma)$ stands for
the multiplicity $\numof{\Set{i|\gamma_i=s}}$ 
of $s\in \set{1,\dots,d}$ in $\gamma$.
\end{proposition}

\begin{remark}
By interpreting the definition of $L_{n,d}$,
we obtain
\begin{align*}
&N(n,d;\gamma)
=\\ &
\Set{
\nu\colon 2^{\Set{1,\ldots,l}}\to\NN
|
\begin{array}{cc}
\displaystyle
\sum_{I\in 2^{\Set{1,\ldots,l}}}\nu(I)=n.
&
\\
\displaystyle
\sum_{I\colon i\in I} \nu(I) = d+1-\gamma_i 
&(\forall i).
\\
\displaystyle
\sum_{J\colon I\subset J} \nu(J)  <  d+1 - \sum_{i\in I}\gamma_i
& (\forall I,\,\numof{I}\ge 2).
\end{array}
}.
\end{align*}
Although 
\eqref{eq:number-type}
is an explicit formula for $\lambda_{n,d}(\gamma)$, 
it is actually difficult to perform the summation in
\eqref{eq:number-type} due to the large size of $N(n,d;\gamma)$.
For example, for $d=6$ and $\gamma=(1,1,1,1,1,1)$.
$\numof{N(n,6;(1,1,1,1,1,1))}=109719496370$. 
\end{remark}

Let
$\tilde N'(j,d;\gamma)=\Set{\tau_T \in \tilde N(j,d;\gamma)|\numof{\tau_T(\emptyset)}=0}$ and
$c(j,d;\gamma)
=\numof{\tilde N'(j,d;\gamma)}$.
Since $\numof{\tau_T(\emptyset)}$
is the number of `symbols' which do not appear in $T$,
$\tilde N'(j,d;\gamma)$
can be identified with $\Set{ T \in \tilde L(j,d;\gamma)| \bigcup_{i} T_i =\Set{1,\ldots,j}}$.
Let $M_{j,n}$ be the set of maps $\sigma\colon \Set{1,\ldots,j} \to \Set{1,\ldots,n}$ 
such that $\sigma(1)<\sigma(2)<\cdots<\sigma(j)$.
Then there exists a natural bijection between 
$\tilde N'(j,d;\gamma) \times M_{j,n}$
and $\tilde N(n,d;\gamma)$.
Hence
\begin{align*}
\numof{\tilde  N(n,d;\gamma)}&=\binom{n}{j}\numof{\tilde N'(j,d;\gamma)}
=\binom{n}{j}c(j,d;\gamma).
\end{align*}
Note that $c(j,d;\gamma)$ does not depend on $n$.
This  equation implies the following proposition.
\begin{proposition}
\label{prop:number-type-modified}
For $\gamma\in\Gamma(n,d)\setminus\Set{\emptypartition}$,
$\lambda_{n,d}(\gamma)$ satisfies 
\begin{align*}
\lambda_{n,d}(\gamma)
&
=\frac{1}{\prod_{s=1}^d m_s(\gamma) !}\sum_{j} c(j,d;\gamma) \binom{n}{j}.
\end{align*}
\end{proposition}
We can also describe $c(j,d;\gamma)$ as the sum of multinomial coefficients:
\begin{align*}
c(j,d;\gamma)
 =
\sum_{\nu' \in N(j,d;\gamma)\colon \nu'(\emptyset)=0} \binom{j}{\nu'(I_1),\ldots,\nu'(I_{2^l-1})},
\end{align*}
where 
$2^{\set{1,\ldots,l}}\setminus \Set{\emptyset}=\Set{I_1,\ldots,I_{2^l-1}}$.

\section{Examples of computations}
\label{sec:tables}

In this section, 
we illustrate the computation of the characteristic polynomial
$\chi_{n,d}(t)$.
We compute the values of $\lambda_{n,d}(\gamma)$ and $c(j,d;\gamma)$
for some $\gamma=(\gamma_1,\ldots,\gamma_l)$ and $d$.
If the length $l$ of
$\gamma=(\gamma_1,\ldots,\gamma_l)$ is small,
then the size of $2^{\{1,\ldots,l\}}$ is small,
so that it is easy to compute $\lambda_{n,d}(\gamma)$ and $c(j,d;\gamma)$.

\subsection{The case when $l=0$}
In this case $\gamma=\emptypartition$  
and by definition 
\[
\lambda_{n,d}(\emptypartition)=\numof{\Set{\minelm}}=1, \quad
\mu_{n,d}(\emptypartition)=\mu_{n,d}(\minelm)=1.
\]
Therefore 
$\mu_{n,0}(\maxelm)=-1$ can also be verified as
\begin{align*}
\mu_{n,0}(\maxelm)
&=-\sum_{i=0}^{0}\sum_{\gamma\partitionof i}\lambda_{n,0}(\gamma)\mu_{n,0}(\gamma)
=-\lambda_{n,0}(\emptypartition)\mu_{n,0}(\emptypartition)=-1.
\end{align*}
This equation implies
\begin{align*}
\mu_{n,d}((\underbrace{1,\ldots,1}_{m}))=
\prod_{i=0}^{m} \mu_{n-d,0}(\maxelm)=(-1)^m.
\end{align*}
We also have
\begin{align*}
\chi(\AAA_{n,0},t)
=
\sum_{i=0}^0 \sum_{\gamma\partitionof i}
\lambda_{n,0}(\gamma)\mu_{n,0}(\gamma)t^{0-i}
=\lambda_{n,0}(\emptypartition)\mu_{n,0}(\emptypartition)=1.
\end{align*}

\subsection{The case when $l=1$}
In this case $\gamma=(\gamma_1)$. Then
$c(j,d;(\gamma_1))=\delta_{j,d+1-\gamma_1}$ and
\[
\lambda_{n,d}((\gamma_1))=\binom{n}{d+1-\gamma_1}.
\]
For example, $\lambda_{n,1}((1))=\binom{n}{1+1-1}=n$.
Hence we have
\begin{align*}
\mu_{n,1}(\maxelm)
&=-\sum_{i=0}^{1}\sum_{\gamma\partitionof i} \lambda_{n,1}(\gamma)\mu_{n,1}(\gamma)\\
&=-\lambda_{n,1}(\emptypartition)\mu_{n,1}(\emptypartition)
-\lambda_{n,1}((1))\mu_{n,1}((1))\\
&=-1\cdot 1 -n\cdot (-1)^1\\
&=n-1.
\end{align*}
This implies
\begin{align*}
\mu_{n,d}((\underbrace{2,\ldots,2}_{m}))=
\prod_{i=0}^{m} \mu_{n+1-d,1}(\maxelm)=(n-d)^m.
\end{align*}
Hence
\begin{align*}
\mu_{n,d}((\underbrace{2,\ldots,2}_{m_2},\underbrace{1,\ldots,1}_{m_1}))
&=\mu_{n,d}((\underbrace{2,\ldots,2}_{m_2}))\mu_{n,d}((\underbrace{1,\ldots,1}_{m_1}))\\
&=(-1)^{m_1}(n-d)^{m_2}.
\end{align*}
We also have
\begin{align*}
\chi(\AAA_{n,1},t)
&=
\lambda_{n,1}((1))\mu_{n,1}((1))+
\lambda_{n,1}(\emptypartition)\mu_{n,1}(\emptypartition)t
=-n+t.
\end{align*}

\subsection{The case when $l=2$}
In this case $\gamma=(\gamma_1,\gamma_2)$.
Let $\nu \in N(j,d;\gamma)$ with $\nu(\emptyset)=0$ and $\nu(\set{1,2})=v$.
Then $\nu(\set{1})=d+1-\gamma_1-v$ and $\nu(\set{2})=d+1-\gamma_2-v$.
Since $\nu(\emptyset)=0$,
$j=\nu(\set{1})+\nu(\set{2})+\nu(\set{1,2})=2d+2-\gamma_1-\gamma_2-v$.
Hence $\nu(\set{1,2})=v=2d+2-\gamma_1-\gamma_2-j$,
$\nu(\set{1})=d+1-\gamma_1-v=\gamma_2+j-d-1$, and
$\nu(\set{2})=d+1-\gamma_2-v=\gamma_1+j-d-1$.
Since $\nu$ satisfies $\nu(\set{1,2})<d+1-\gamma_1-\gamma_2$,
we also have $j>d+1$, which implies $\nu(\set{1}),\nu(\set{2})\geq 0$.
Since $\gamma(\set{1,2})$ should be nonnegative,
we have $2d+2-\gamma_1-\gamma_2\geq j$.
Hence 
\begin{align*}
c(j,d;(\gamma_1,\gamma_2))&=
\binom{j}{\gamma_1+j-d-1,\gamma_2+j-d-1,2d+2-\gamma_1-\gamma_2-j}
\end{align*}
if $ 2d+2-\gamma_1-\gamma_2\geq j>d+1$.
Let $j'=j-d-1$. Then we can rewrite the equation as
\begin{align*}
c(d+1+j',d;(\gamma_1,\gamma_2))
&=
\binom{d+1+j'}{\gamma_1+j',\gamma_2+j',d+1-\gamma_1-\gamma_2-j'}\\
&=
\binom{d+1+j'}{\gamma_1+j'}\binom{d+1-\gamma_1}{\gamma_2+j'}
\end{align*}
if $ d+1-\gamma_1-\gamma_2\geq j'>0$.
Hence 
\begin{align*}
\lambda_{n,d}((\gamma_1, \gamma_2)) =
\sum_{j'=1}^{d+1-\gamma_1-\gamma_2}
\binom{d+1+j'}{\gamma_1+j'}\binom{d+1-\gamma_1}{\gamma_2+j'}
\binom{n}{d+1+j'}
\end{align*}
if $\gamma_1\neq \gamma_2$.
We can also obtain
\begin{align*}
\lambda_{n,d}((\gamma_1, \gamma_2)) =
\frac{1}{2}\sum_{j'=1}^{d+1-\gamma_1-\gamma_2}
\binom{d+1+j'}{\gamma_1+j'}\binom{d+1-\gamma_1}{\gamma_2+j'}
\binom{n}{d+1+j'}
\end{align*}
if $\gamma_1= \gamma_2$.
For example, 
\begin{align*}
\lambda_{n,2}((1,1)) 
&=
\frac{1}{2}\sum_{j'=1}^{2+1-2} \binom{3+j'}{1+j'}\binom{3-1}{1+j'}\binom{n}{3+j'}
=3\binom{n}{4}.
\end{align*}
We also have
\begin{align*}
\mu_{n,2}(\maxelm)
=&-\lambda_{n,2}(\emptypartition)\mu_{n,2}(\emptypartition)
-\lambda_{n,2}((1))\mu_{n,2}((1))
\\&
-\lambda_{n,2}((2))\mu_{n,2}((2))
-\lambda_{n,2}((1,1))\mu_{n,2}((1,1))
\\
=&-1\cdot 1
-\binom{n}{2}\cdot(-1)^1
-\binom{n}{1}\cdot (n-2)^1
-3\binom{n}{4}\cdot (-1)^2  
\\
=&-1+\frac{-n(n-3)}{2}-3\binom{n}{4}
=
-1+n-\binom{n}{2}-3\binom{n}{4}\\
=&\frac{-(n- 2)(n- 1)(n^2 - 3n + 4)}{8}.
\end{align*}
This implies
\begin{align*}
\mu_{n,d}((\underbrace{3,\ldots,3}_{m}))
&=\prod_{i=0}^{m} \mu_{n+2-d,2}(\maxelm)\\
&=\left(\frac{-(n - d)(n - d + 1)(d^2 - 2dn + n^2+ n + 2- d)}{8}\right)^m.  
\end{align*}
If $1\leq \gamma_i\leq 3$ for all $i$,
we can obtain an explicit formula for 
$\mu_{n,d}(\gamma)$.
We also have
\begin{align*}
\chi(\AAA_{n,2},t)
=&
\lambda_{n,2}((1,1))\mu_{n,2}((1,1))+
\lambda_{n,2}((2))\mu_{n,2}((2))
\\&
+\lambda_{n,2}((1))\mu_{n,2}((1))t+
\lambda_{n,2}(\emptypartition)\mu_{n,2}(\emptypartition)t^2
\\
=&
3\binom{n}{4}
+(n-2)n
-\binom{n}{2}t+
t^2\\
=&
-n
+2\binom{n}{2}
+3\binom{n}{4}
-\binom{n}{2}t+t^2.
\end{align*}

\subsection{The case when $l\geq 3$}
Finally 
we consider the case when $l\geq 3$.
In this case, it seems hard to compute $c(j,d;\gamma)$
by hand,
because $\tilde N'(j,d;\gamma)$ is complicated 
due to the size of $2^{\set{1,\ldots,l}}$.
We can, however, enumerate $c(j,d;\gamma)$
for a particular $j$, $d$ and $\gamma$
by a computer.
For example, Program \ref{prog:sage} is a program for 
the computer algebra system \path{Sage}\cite{sage}
to define the function \path{getC} to compute
the value of $c(j,d;\gamma)$.

\begin{program}%
\label{prog:sage}
\footnotesize
\begin{verbatim}
def condition2(n,d,gamma,nu):
  l=len(gamma)
  for i in range(l):
    lhs=sum(nu[I] for I in nu.keys() if i in I)
    rhs=d+1-gamma[i]
    if not lhs == rhs:
      return False
  return True
    
def condition3(n,d,gamma,nu):
  for I in nu.keys():
    if len(I)>1:
      rhs=sum(nu[J] for J in nu.keys() if I.issubset(J))
      lhs=d+1-sum(gamma[i] for i in I)
      if not rhs < lhs:
        return False
  return True

def xLatticePoints(n,total):
  if n==1:
    yield [total]
    return
  for i in (0..total):
    for v in xLatticePoints(n-1,total-i):
      yield [i]+v

def xN(j,d,gamma):
  l=len(gamma)
  index=[frozenset(I) for I in powerset((0..l-1)) if I != []]
  for v in xLatticePoints(2^l-1,j):
    nu=dict(zip(index,v))
    if condition2(j,d,gamma,nu):
      if condition3(j,d,gamma,nu):
        yield nu

def getC(j,d,gamma):
  return sum(multinomial(nu.values()) for nu in xN(j,d,gamma))
\end{verbatim}
\end{program}
\normalsize

Running the program,
we can obtain the table of $c(j,d;\gamma)$.
However computation is heavy 
if $j$, $d$ or the length of $\gamma$ is large.
For example, the following shows 
the values of $c(3d,d;(1,1,1))$ and its processing time:
\begin{align*}
c(9,3;(1,1,1))&= 1680 &0.44S&\\
c(12,4;(1,1,1))&= 34650 &1.36S&\\
c(15,5;(1,1,1))&= 756756 &3.50S&\\
c(18,6;(1,1,1))&= 17153136 &7.87S&\\
c(21,7;(1,1,1))&= 399072960 &15.96S&\\
c(24,8;(1,1,1))&= 9465511770 &30.04S&\\
c(27,9;(1,1,1))&= 227873431500 &53.10S&.
\end{align*}
The length of $\gamma$ is more critical than $j$ and $d$.
For example,
we enumerated $c(8,4;(1,1,1))=68040$ in $0.3$ seconds and
$c(8,4;(1,1,1,1))=2900520$ in $1$ minute.  
Other computational data by the program are available 
from the authors.

\section{Combinatorics of discriminantal arrangements}
\label{sec:disc} 
We now consider the arrangement $\BBB(d+k+1,k)=\AAA_{d+k+1,d}$,
parametrized by $d$ for fixed $k$.
The combinatorics is already very hard for $k=2$. Therefore 
in this section we discuss the class of arrangements $\AAA_{d+k+1,d}$ for $k=0,1,2$.

We first show an inequality for general  $k\geq 0$.
Each element $S \in L^{d+k+1,d}$ satisfies
the condition
\begin{align*}
D^k(S')>0
\end{align*}
for $S'\subset S$ with $\numof{S'}>1$.
For $S_i,S_j \in S$ with $S_i\neq S_j$,
the condition means
\begin{align*}
D^k(\Set{S_i,S_j})
&=\numof{S_i \cup S_j} -k -(\numof{S_i}-k) -(\numof{S_j}-k)\\
 &=\numof{S_i \cup S_j} -\numof{S_i}-\numof{S_j}+k\\
 &=-\numof{S_i \cap S_j}+k>0.
\end{align*}
Hence we have the inequality
\begin{align}
\label{ineq:pair}
k>\numof{S_i \cap S_j}\geq 0.
\end{align}

\subsection{The case when $k=0$}
Here we consider $L^{d+1,d}$.
We show that
$L^{d+1,d}$
is isomorphic to the intersection lattice of
the Boolean arrangement, i.e., 
the set of the hyperplanes
$\Set{(x_1,\ldots,x_{d+1})|x_i=0}$ for $i=1,\ldots,d+1$,
as posets.

For $k=0$ the inequality \eqref{ineq:pair} is
\begin{align*}
0>\numof{S_i \cap S_j}\geq 0,
\end{align*}
which is always false.
Hence $\numof{S}\leq 1$ for  $S \in L^{d+1,d}$.
Therefore we have
\begin{align*}
L^{d+1,d}
&=\Set{\{S_1\} | S_1\subset [d+1], \numof{S_1}>0} \cup \Set{\emptyset}.
\end{align*}
The type of any element 
in  $L^{d+1,d}$
has at most one row,
namely,
\begin{align*}
\Gamma(d+1,d)=\Set{\emptypartition,(1),(2),\ldots,(d+1)}.
\end{align*}
For $\Set{S_1}$ and $\Set{S'_1} \in L^{d+1,d}$,
$\Set{S_1} \leq \Set{S'_1}  $
if and only if $S_1 \subset S'_1$.
Hence the map from $L^{d+1,d}$ to 
the Boolean lattice $2^{[d+1]}$ defined by
\begin{align*}
\{S_1\}   &\mapsto S_1,\\
\emptyset &\mapsto \emptyset
\end{align*}
is an order-preserving bijection.
Via this bijection,
 $L^{d+1,d}$ is isomorphic to
the Boolean lattice $2^{[d+1]}$
as posets.

\subsection{The case when $k=1$}
Here we consider $L^{d+2,d}$.
It is known that 
$L^{d+2,d}$ is isomorphic 
to intersection poset of the braid arrangement (e.g.\ Section 1 of
Falk\cite{falk1994}, Section 5.6 of Orlik and Terao\cite{orlik-terao-book}).
We summarize the correspondence here.

In this case, \eqref{ineq:pair} implies
$\numof{S_i \cap S_j} = 0$ for 
$S_i, S_j \in S \in L^{d+2,d}$ with $S_i \neq S_j$.
Hence 
\begin{align*}
 L^{d+2,d}=
\Set{
 S \subset 2^{[d+2]} | 
 \begin{array}{c}
 \text{$S_i\cap S_j = \emptyset$ for $S_i,S_j \in S$ with $S_i\neq S_j$.}\\
 \text{$1<\numof{S_i}<d+2$ for $S_i \in S$.} 
 \end{array}
}.
\end{align*}

Since $S_i \cap S_j =\emptyset$ for $S_i\neq S_j\in S$,
we have $\sum_{S_i\in S} \numof{S_i}=\numof{\bigcup_{S_i\in S} S_i}$ 
for $S\in L^{d+2,d}$.
Hence,
for $S \in L^{d+2,d}$ of type $(\gamma_1,\ldots,\gamma_l)$,
we have
\begin{align*}
 \sum_i \gamma_i
 = \sum_{S_i  \in S} \codim^1 (S_i)
 = \sum_{S_i  \in S} (\numof{S_i}-1)
 = -l+\numof{\bigcup_{S_i  \in S} S_i}.
\end{align*}
Since $\bigcup_{S_i  \in S} S_i \subset [d+2]$, we obtain
\begin{align*}
 \sum_i \gamma_i = -l+\numof{\bigcup_{S_i  \in S} S_i} \leq -l+d+2.
\end{align*}
Therefore,
 if $\sum_{i=1}^{l} \gamma_i > d+2-l$,
then there does not exist
 an element of $L^{d+2,d}$ of type $(\gamma_1,\ldots,\gamma_l)$.
 Moreover, if $\sum_i \gamma_i \leq d+2-l$,
 then the number of 
 elements of $L^{d+2,d}$
 whose type is $(\gamma_1,\ldots,\gamma_l)$
 is
 \begin{align*}
 \frac{1}{\prod_s m_s(\gamma)!}\binom{d+2}{\gamma_1,\ldots,\gamma_l,d+2-\sum_i \gamma_i }.
 \end{align*}

Now we discuss the structure of the poset $L^{d+2,d}$.
Let $H_{j,j'}=\Set{(x_1,\ldots,x_{d+2}) \in \RR^{d+2}| x_j = x_{j'}}$.
The arrangement $\Set{H_{j,j'}| j\neq j'}$ is called the braid arrangement
$\AAA_{d+1}$.
For $S \in L^{d+2,d}$,
define 
\begin{align*}
\varphi(S)= \Set{(x_1,\ldots,x_{d+2}) \in \RR^{d+2}| \text{$x_j=x_{j'}$ for all $j,j' \in S_i \in S$}}.
\end{align*}
If $S=\emptyset$, then $\varphi(S)=\RR^{d+2}$.
If $S=\Set{S_1}$ and $\numof{S_1}\geq 2$, then 
\begin{align*}
\varphi(S)&=\Set{(x_1,\ldots,x_{d+2}) \in \RR^{d+2}| x_j=x_{j'}\ (j,j' \in S_1)}\\
&=\bigcap_{j, j' \in S_1} H_{j,j'} \in \AAA_{d+1}.
\end{align*}
If $S=\Set{S_1,\ldots,S_l}$, then 
\begin{align*}
\varphi(S)=\bigcap_{S_i \in S} \varphi(\Set{S_i}) \in \AAA_{d+1}
\end{align*}
since $S_i \cap S_{i'} =\emptyset$ for $S_i \neq S_{i'} \in S$.
Hence $\varphi$ is a map from $L^{d+2,d}$ to the
intersection lattice of the braid arrangement $\AAA_{d+1}$.
Conversely, for $H_{i_1,j_1}\cap \cdots \cap H_{i_m,j_m}$,
there exist a set partition $S=\Set{S_1,\ldots,S_l}$ of $[d+2]$
such that 
\begin{align*}
H_{i_1,j_1}\cap \cdots \cap H_{i_m,j_m}
= \Set{(x_1,\ldots,x_{d+2}) \in \RR^{d+2}| x_j=x_{j'} \ (j,j' \in S_i \in S)}.
\end{align*}
Let 
$\psi(H_{i_1,j_1}\cap \cdots \cap H_{i_m,j_m})=\Set{S_i \in S|\numof{S_i}\geq 2}$.
Then,
since $S_i\cap S_{i'}=\emptyset$ for $S_i\neq S_{i'} \in S$
and $2 \leq \numof{S_i} \leq d+2$ for $S_i \in \psi(H_{i_1,j_1}\cap \cdots \cap H_{i_m,j_m})$,
 $\psi(H_{i_1,j_1}\cap \cdots \cap H_{i_m,j_m})$ is an element in 
 $L^{d+2,d}$.
Since the map $\psi$ provides the inverse map of $\varphi$,
the map $\varphi$ is a bijection.
For $S_1\subset S'_1 \subset [d+2]$ with $\numof{S_1}>1$,
we have
\begin{align*}
\varphi(\Set{S_1})=\bigcap_{j,j' \in S_1} H_{j,j'}
\supset 
\bigcap_{j,j' \in S_1} H_{j,j'} \cap \bigcap_{j,j' \in S'_1 \setminus S_1} H_{j,j'}
=\varphi(\Set{S'_1}).
\end{align*}
Hence, for 
$S=\Set{S_1,\ldots,S_l} < S'= \Set{ S'_1,\ldots,S'_{l'}} \in L^{d+2,d}$,
we have
\begin{align*}
\varphi(S)=\bigcap_{S_i \in S} \varphi(\Set{S_i}) \supset
\bigcap_{S'_i \in S'} \varphi(\Set{S'_i})= 
\varphi(S'),
\end{align*}
and $\varphi$ is an order-preserving bijection.
Therefore $L^{d+2,d}$ is isomorphic to 
the intersection lattice of the braid arrangement $\AAA_{d+1}$
as posets.

\subsection{The case when $k=2$.}
Here we consider $L^{d+3,d}$.
Although $L^{d+3,d}$ is an analogue of 
the Boolean arrangement and the braid arrangement, 
the structure of $L^{d+3,d}$
is complicated and it is difficult to compute  the characteristic polynomials
for large $d$.
For $S=\Set{S_1,\ldots,S_l} \in L^{d+k+1,d}$,
the set
\begin{align*}
\Delta(S)=\Set{ I \subset [l] | \bigcap_{i \in I} S_i \neq \emptyset}
\end{align*}
is a simplicial complex,
or equivalently, $\Delta(S)$ satisfies
\begin{align*}
F\in\Delta(S),\ F' \subset F \implies F'\in\Delta(S).
\end{align*}
For  $k=2$
the inequality \eqref{ineq:pair} 
implies $\numof{S_i \cap S_j} \leq 1$.
Hence any two facets  of the simplicial complex $\Delta(S)$
share at most one vertex.
Moreover each face in the simplicial complex with positive dimension
corresponds to only one symbol in $[d+3]$.
Here we discuss a way to enumerate $L^{d+3,d}$ 
from the information of the
simplicial complex $\Delta(S)$.

First we define notation for simplicial complexes.
Let $\VVV_0(l)$ be the set of simplicial complexes $\Delta\subset 2^{[l]}$ 
such that $\Set{i} \in \Delta$ for all $i$.
For $\Delta \in \VVV_0(l)$,
we define $\FFF(\Delta)$ to be the set 
of facets with positive dimension.
Let $\VVV_1(l)$ be the set of  simplicial complexes $\Delta\subset 2^{[l]}$ 
 such that 
the dimension of the face $F\cap F'$ is at most one
for any distinct facets $F,F' \in \FFF(\Delta)$.
For $I\subset [l]$ and $\Delta \in \VVV_1(l)$,
define $\Delta|_{I} = \Set{ F\cap I| F \in \Delta}$.
Let 
\begin{align*}
\VVV_2(l)=\Set{
\Delta \in \VVV_1(l)
|
2(\numof{I}-1)>\sum_{F\in \FFF(\Delta|_{I})} \dim F
\quad(\forall I\subset [l], \numof{I}>1).
}.
\end{align*}
For $\Delta \in \VVV_2(l)$,
we define $\alpha(\Delta)=(\alpha_1,\ldots,\alpha_l)$
with $\alpha_i = \numof{\Set{F \in \FFF(\Delta)| i \in F}}$.
\begin{example}
For $\Delta \in \VVV_2(2)$,
$\FFF(\Delta)=\Set{\{1,2\}}$ or $\emptyset$.
\end{example}

\begin{example}
For $\Delta \in \VVV_2(3)$,
$\FFF(\Delta)$
is a subset of one of the following:
\begin{align*}
\Set{\{1,2,3\}}, \Set{\{1,2\},\{2,3\},\{1,3\}}.
\end{align*}
Figure~\ref{fig:2:1} and \ref{fig:2:2}
show $\Set{\{1,2,3\}}$ and $\Set{\{1,2\},\{2,3\},\{1,3\}}$
as hypergraphs, respectively, where a hyperedge $e$ with $\numof{e}=2$ is
written as an ordinary edge.
\begin{figure}[t]
\hfil
\subfigure[]{\scIIIii\label{fig:2:1}}
\hfil
\subfigure[]{\scIIIi\label{fig:2:2}}
\caption{Example of $\FFF(\Delta)$ for $\Delta\in \VVV_2(3)$}
\label{fig:exampleFDelta:3}
\end{figure}
\end{example}

\begin{example}
For $\Delta \in \VVV_2(4)$,
$\FFF(\Delta)$
is a subset of one of the following:
\begin{align*}
&\Set{\{\sigma_1,\sigma_2,\sigma_3,\sigma_4\}} 
&&(\text{Figure \ref{fig:3:1}}), \\
&\Set{\{\sigma_1,\sigma_2,\sigma_3\},\{\sigma_1,\sigma_4\},\{\sigma_2,\sigma_4\},\{\sigma_3,\sigma_4\}}
&&(\text{Figure \ref{fig:3:2}}),\\
&\Set{\{\sigma_1,\sigma_2\},\{\sigma_2,\sigma_3\},\{\sigma_1,\sigma_3\},\{\sigma_4,\sigma_2\},\{\sigma_2,\sigma_3\},\{\sigma_4,\sigma_3\}}
&&(\text{Figure \ref{fig:3:3}}),
\end{align*}
where $\sigma \in \mathfrak{S}_4$.
\begin{figure}[t]
\hfil
\subfigure[]{\scIVii\label{fig:3:1}}
\hfil
\subfigure[]{\scIViii\label{fig:3:2}}
\hfil
\subfigure[]{\scIVi\label{fig:3:3}}
\caption{Example of $\FFF(\Delta)$ for $\Delta\in \VVV_2(4)$}
\label{fig:exampleFDelta:4}
\end{figure}
\end{example}
\begin{example}
For $\Delta \in \VVV_2(5)$,
the hyper graph $\FFF(\Delta)$
is a subset of one of (a)--(i) 
in Figure \ref{fig:exampleFDelta:5}.
\begin{figure}[t]
\hfil
\subfigure[]{\scVi}
\hfil
\subfigure[]{\scVii}
\hfil
\subfigure[]{\scViii}
\hfil
\subfigure[]{\scViv}\\
\hfil
\subfigure[]{\scVv}
\hfil
\subfigure[]{\scVvi}
\hfil
\subfigure[]{\scVvii}
\hfil
\subfigure[]{\scVviii}
\hfil
\subfigure[]{\scVix}
\caption{Example of $\FFF(\Delta)$ for $\Delta\in \VVV_2(5)$}
\label{fig:exampleFDelta:5}
\end{figure}
\end{example}

\newcommand{\restricteddelta}{\Delta{|}_{[l-1]}}
\begin{remark}
It follows
by definition
that
\begin{align*}
\VVV_2(l-1)=
\left\{ \restricteddelta \bigm|  \Delta \in \VVV_2(l) \right\}.
\end{align*}
\end{remark}

\begin{remark}
\label{rem:graphcasemax}
Consider a graph $\Delta$ in $\VVV_2(l)$, namely, $\Delta\in \VVV_2(l)$ such that 
$\dim F =1$ for $F\in \FFF(\Delta)$.
Since $\sum_{F\in \FFF(\Delta)} \dim F = \numof{\FFF(\Delta)}$,
the number of edges in $\Delta$ is 
less than or equal to $2l-3$.
Let $\Delta$ be the graph with edges
\begin{align*}
\begin{cases}
\Set{1,j} &(j=2,\ldots,l)\\
\Set{j,j+1} &(j=2,\ldots,l-1).
\end{cases}
\end{align*}
See Figure \ref{fig5}. Then $\numof{\Delta}=2l-3$.
If a subset $I\subset [l]$ does not contain $1$,
then $\FFF(\Delta|_I)=\Set{\{i,i+1\}| i,i+1 \in I}$,
which implies $\numof{\FFF(\Delta|_I)} \leq \numof{I}-1$. 
If a subset $I\subset [l]$ contains $1$,
then $\FFF(\Delta|_I)=\Set{\{i,i+1\}| i,i+1 \in I\setminus\Set{1}}\cup \Set{\{i,1\}| i \in I\setminus\Set{1}}$,
which implies $\numof{\FFF(\Delta|_I)} \leq (\numof{I}-2)+(\numof{I}-1)=2\numof{I}-3$.
Hence $\Delta \in \VVV_2(l)$.
Therefore
the inequality $\numof{\FFF(\Delta)}\leq 2l-3$ for $\Delta\in\VVV_2(l)$
is the best possible.
\begin{figure}[t]
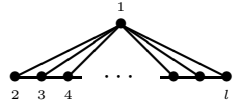

\hfil\figGraph
\caption{$\Delta$ in Remark \ref{rem:graphcasemax}}
\label{fig5}
\end{figure}
\end{remark}

\begin{remark}
Let $\Delta\in\VVV_2(l)$ and $F=\Set{v_1,\ldots,v_m}\in\FFF(\Delta)$.
Consider the simplicial complex $\Delta'$ such that
\begin{align*}
\FFF(\Delta')=(\FFF(\Delta)\setminus F)\cup \Set{\{v_i,v_{i+1}\}|i=1,\ldots,m}.
\end{align*}
Since $\dim F|_I\geq \sum_{v_i,v_{i+1}\in I} \dim(\{v_i,v_{i+1}\})$, 
the simplicial complex $\Delta'$ is in $\VVV_2(l)$.
It also follows that $\numof{\FFF(\Delta)}\leq \numof{\FFF(\Delta')}$.
Hence there exists a graph $\Delta\in\VVV_2(l)$
such that $\numof{\FFF(\Delta)}=\max\Set{\numof{\FFF(\Delta')}|\Delta'\in\VVV_2(l)}$.
Therefore we have
\begin{align*}
\max\Set{\numof{\FFF(\Delta')}|\Delta'\in\VVV_2(l)}=2l-3 .
\end{align*}
\end{remark}

Now we continue the investigation of $L^{d+3,d}$.
For the computational reason,
we consider the tuples of subsets in $[d+3]$.
Let 
\begin{align*}
&\tilde L^{d+3,d}=\Set{(S_1,\ldots,S_l)| \Set{S_1,\ldots,S_l}\in L^{d+3,d}}, \\
&\tilde L^{d+3,d;\gamma}=\Set{(S_1,\ldots,S_l)\in
\tilde L^{d+3,d}| \numof{S_i}=\gamma_i+2}.
\end{align*}
For  $S=(S_1,\ldots,S_l) \in \tilde L^{d+3,d}$ and $t \in [d+3]$,
we define 
\begin{align*}
X(S,t) &=\Set{ i | t \in S_i }.
\end{align*}
For each facet $F$ of $\Delta(S)$, 
there exists a symbol $t \in [d+3]$ such that $F=X(S,t)$.
The inequality \eqref{ineq:pair} implies
\begin{align*}
2>\numof{S_i \cap S_j} \geq 0
\end{align*}
for $S_i\neq S_j \in S \in L^{d+3,d}$.
Hence 
$\numof{S_i \cap S_j}$ is at most one,
which implies
$\numof{X(S,t)\cap X(S,t')}\leq 1$ for $t\neq t'$.
In other words,
the dimension of the face
$X(S,t) \cap X(S,t')$ of $\Delta(S)$
is at most one.
Hence we have $\Delta(S)\in \VVV_1(l)$.
Moreover, there exists one-to-one correspondence between 
$\FFF(\Delta(S))$
and $\Set{t \in [d+3] | \numof{X(S,t)}>1}$.

The number of symbols used in $S\in \tilde L^{d+3,d}$
is $\numof{\bigcup_i S_i}$.
On the other hand,
since 
the number of times the symbol $t$ appears in $S$ is $\numof{X(S,t)}$,
the number of symbols used in $S$ is
$\sum_i \numof{S_i} - \sum_t(\numof{X(S,t)}-1)$.
Hence we have 
$\numof{\bigcup_i S_i}=\sum_i \numof{S_i} - \sum_t(\numof{X(S,t)}-1)$.
Therefore for any $k\ge 0$
\begin{align*}
D^{k}(\Set{S_1,\ldots,S_l})
&=\Bigr|\bigcup_{i \in [l]} S_i \Bigr| -k - \sum_{i\in [l]} (\numof{S_i}-k)\\
&=\sum_{i\in [l]} \numof{S_i} - \sum_{t\in [d+3]}(\numof{X(S,t)}-1)
-k - \sum_{i\in [l]} (\numof{S_i}-k)\\
&=k(l-1) - \sum_t(\numof{X(S,t)}-1).
\end{align*}
Hence 
the inequality 
$D^{k}(\Set{S_1,\ldots,S_l})>0$ for $k=2$ 
is rewritten as
\begin{align*}
2(l-1)>\sum_{t \in [d+3]}\dim X(S,t),
\end{align*}
where  $\dim X(S,t)$ 
stands for the dimension $\numof{X(S,t)}-1$ of $X(S,t)$ as a simplex.
Since there exists a one-to-one correspondence 
between $\FFF(\Delta(S))$
and $\Set{t \in [d+3] | \dim X(S,t)>0}$,
the inequality means
\begin{align*}
2(l-1)>\sum_{F\in \FFF(\Delta(S))} \dim F.
\end{align*}
Similarly we have inequalities 
\begin{align*}
2(\numof{I}-1)>\sum_{F\in \FFF(\Delta(S)|_{I})} \dim F
\end{align*}
for $I \subset [l]$ with $\numof{I}>1$.
Hence 
the simplicial complex $\Delta(S)$ is in $\VVV_2(l)$.

Conversely, let $\Delta \in \VVV_2(l)$ with $\alpha(\Delta)=(\alpha_1,\ldots,\alpha_l)$ and $S_1,\ldots,S_l \subset [d+3]$.
Assume that
\begin{align*}
&S_i \cap S_j = \emptyset  \quad  (i\neq j) , \\
&\numof{\FFF(\Delta)}+\sum_{i=1}^{l} \numof{S_i} \leq d+3 ,\\
&3\leq  \alpha_i+\numof{S_i} \leq d+3.
\end{align*}
Fix an injection $\tau$  from $\FFF(\Delta)$ to $[d+3]\setminus \bigcup_i S_i$,
and let $S^{\tau}_i = \Set{\tau(F) | i\in F\in \FFF(\Delta)} \cup S_i$.
From these assumptions we have 
$3\leq  \numof{S^{\tau}_i} \leq d+3$.
Moreover, since $\Delta \in \VVV_2(l)$,
we have 
$D^{2}(\Set{S^{\tau}_{i_1},\ldots,S^{\tau}_{i_{l'}}})>0$.
Hence $(S_1,\ldots,S_l) \in  \tilde L^{d+3,d}$.
By this construction, we can enumerate the number of elements in 
$\tilde L^{d+3,d}$ of type $\gamma$ as follows.

\begin{proposition}
\label{prop:discriminantal-type}
For $\gamma=(\gamma_1,\ldots,\gamma_l)$,
\begin{align}
&\numof{\tilde L^{d+3,d;\gamma}}
\nonumber\\
& =
\sum_{\alpha'}
\sum_{\Delta }
\binom{d+3}{\alpha'_1,\ldots,\alpha'_l, d+3-\sum_i\alpha'_i}
\frac{(d+3-\sum_i\alpha'_i) !}{ (d+3-\sum_i\alpha'_i-\numof{\FFF(\Delta)})!},
\label{eq:discriminantal-type}
\end{align}
where 
the first sum is over $\alpha'=(\alpha'_1,\ldots,\alpha'_l)$
such that
 $0\leq \alpha_i' \leq \gamma_i+2$ for all $i$,
and the second sum is over 
$\Delta \in  \VVV_3(l)$
such that 
 $\alpha(\Delta)= (\gamma_1+2-\alpha_1',\ldots,\gamma_l+2-\alpha_l')$
and $\numof{\FFF(\Delta)}\leq d+3 -\sum_i \alpha'_i$.
\end{proposition}

\begin{remark}
\label{rem:alltypeappear}
In the case when $k\leq 1$, 
there exists $(\gamma_1,\ldots,\gamma_l)$
which does not appear as 
the type of any element  of $L^{d+k+1,d}$.
In the case when $k=2$, 
all $(\gamma_1,\ldots,\gamma_l)\partitionof d' \leq d$
appear as the type of some element  of $L^{d+k+1,d}$.
This can be seen as follows. Consider the simplicial complex $\Delta$
with facets
\begin{align*}
&\Set{ i \in [l]| \text{$i$ is odd} },&
&\Set{ i \in [l]| \text{$i$ is even} },&
&\Set{ j,j+1  } \text{for $j=1,\ldots,l-1$}.
\end{align*}
\begin{figure}[t]
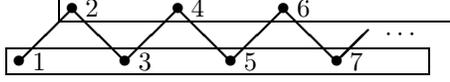

\hfil\scUniv
\caption{Facets of $\Delta$ in Remark \ref{rem:alltypeappear}}\label{fig4}
\end{figure}
See Figure \ref{fig4}. Then 
$\sum_{F\in \FFF(\Delta)} \dim F=2l-3<2(l-1)$
since $\dim \Set{ i \in [l]| \text{$i$ is odd} }+\dim \Set{ i \in [l]| \text{$i$ is even} }=l-2$.
Let $I\subset [l]$ with $\numof{I}\geq 2$.
Then facets of $\Delta|_I$ with positive dimension are
\begin{align*}
&\Set{ i \in I| \text{$i$ is odd} },&
&\Set{ i \in I| \text{$i$ is even} },&
&\Set{ j,j+1 }  \text{for $j,j+1\in I$}.
\end{align*}
Hence 
$\sum_{F\in \FFF(\Delta|_I)} \dim F\leq \numof{I}-2+\numof{I}-1=2\numof{I}-3$.
Therefore $\Delta$ is in $\VVV_2(l)$.
It follows from the direct calculation that
$\alpha_1=\alpha_l=2$ and $\alpha_2=\cdots=\alpha_{l-1}=3$
for $\alpha(\Delta)=(\alpha_1,\ldots,\alpha_l)$.
Since $\numof{\FFF(\Delta)}=2+l-1=l+1$ and $\sum_i \alpha_i =3l-2$,
we have
\begin{align*}
\numof{\FFF(\Delta)}+\sum_{i=1}^{l} (\gamma_i+2-\alpha_i)
=3+\sum_{i=1}^{l}\gamma_i
\end{align*}
for 
$\gamma=(\gamma_1,\ldots,\gamma_l)$.
If  $\sum_i \gamma_i \geq d$,
then we have $\numof{\FFF(\Delta)}+\sum_{i=1}^{l} (\gamma_i+2-\alpha_i) \leq d+3$.
Since $\alpha_i \leq 3$ for each $i$,
we can construct
an element in $L^{d+3,d}$
of type $(\gamma_1,\ldots,\gamma_l)\partitionof d' \leq d$. 
\end{remark}

\section{Computation of $\AAA_{d+3,d}$}
\label{sec:computation:disc}
Here we give some computational data of  
$\Delta \in \VVV_2(l)$ and $\AAA_{d+3,d}$.
More comprehensive computational data are available from the authors.

The summand in the equation \eqref{eq:discriminantal-type}
is written as
\begin{align*}
&\binom{d+3}{\alpha'_1,\ldots,\alpha'_l, d+3-\sum_i\alpha'_i}
\frac{(d+3-\sum_i\alpha'_i) !}{ (d+3-\sum_i\alpha'_i-\numof{\FFF(\Delta)})!}
\\
=&
\binom{d+3}{\alpha'_1,\ldots,\alpha'_l, \numof{\FFF(\Delta)},d+3-\sum_i\alpha'_i-\numof{\FFF(\Delta)}}
\numof{\FFF(\Delta)}!\\
=&
\binom{d+3}{ \numof{\FFF(\Delta)} +\sum_i\alpha'_i}
\binom{ \numof{\FFF(\Delta)}  +\sum_i\alpha'_i}
{\alpha'_1,\ldots,\alpha'_l,\numof{\FFF(\Delta)}}
\numof{\FFF(\Delta)}!.
\end{align*}
For $(\alpha'_1,\ldots,\alpha'_l)$, $(\alpha_1,\ldots,\alpha_l)$ and
$t$, 
define
\begin{align*}
P(\alpha'_1,\ldots,\alpha'_l;t)&=
\binom{d+3}{ t +\sum_i\alpha'_i}
\binom{ t  +\sum_i\alpha'_i}
{\alpha'_1,\ldots,\alpha'_l,t}t!,\\
\VVV_2((\alpha_1,\ldots,\alpha_l),t)
&=\Set{\Delta \in \VVV_2(l)| 
\begin{array}{c}
\alpha(\Delta)=(\alpha_1,\ldots,\alpha_l)\\
|\FFF(\Delta)|=t
\end{array}
}.
\end{align*}
Since $\alpha_i=\gamma_i+2-\alpha'_i$ holds in the equation \eqref{eq:discriminantal-type},
we rewrite the equation
as
\begin{align*}
\numof{\tilde L^{d+3,d;\gamma}}
&=
 \sum_{\alpha}
\sum_t \sum_{\Delta \in \VVV_2(\alpha,t)}
P(\gamma_1+2-\alpha_1,\ldots,\gamma_l+2-\alpha_l;\numof{\FFF(\Delta)})
\\
&=
 \sum_{\alpha}
\sum_t \numof{\VVV_2(\alpha,t)}
P(\gamma_1+2-\alpha_1,\ldots,\gamma_l+2-\alpha_l;t),
\end{align*}
where the sum is over all $\alpha=(\alpha_1,\ldots,\alpha_l)$
such that $0\leq \alpha_i \leq \gamma_i+2$ for all $i$.
Since $\numof{ L^{d+3,d;\gamma}}=\numof{\tilde L^{d+3,d;\gamma}}/\prod_s m_s(\gamma)!$,
by computing $\numof{\VVV_2(\alpha,t)}$,
we obtain 
Table \ref{tab1} of $\lambda_{d+3,d}(\gamma)=\numof{ L^{d+3,d;\gamma}}$.
The table of $\numof{\VVV_2(\alpha,t)}$ is available from the authors.
\begin{table}[htbp]
\tbl{$\lambda_{d+3,d}(\gamma)$}
{\footnotesize
\begin{tabular}{crrrrrrr}
$\gamma$ &$d= 1 $ &$d= 2 $ &$d= 3 $ &$d= 4 $ &$d= 5
$ &$d= 6 $ &$d= 7 $ \\ 
\hline
$ ( 1 )$ & 4 & 10 & 20 & 35 & 56 & 84 & 120
\\ 
\hline
$ ( 2 )$ & 1 & 5 & 15 & 35 & 70 & 126 & 210
\\ 
$ ( 1 ^{ 2 } )$ & 0 & 15 & 100 & 385 & 1120 &
2730 & 5880 \\ 
\hline
$ ( 3 )$ & 0 & 1 & 6 & 21 & 56 & 126 & 252
\\ 
$ ( 2 1 )$ & 0 & 0 & 60 & 455 & 1960 & 6300
& 16800 \\ 
$ ( 1 ^{ 3 } )$ & 0 & 0 & 120 & 1575 & 10080 &
44380 & 154000 \\ 
\hline
$ ( 4 )$ & 0 & 0 & 1 & 7 & 28 & 84 & 210 \\ 
$ ( 3 1 )$ & 0 & 0 & 0 & 105 & 896 & 4284 &
15120 \\ 
$ ( 2 ^{ 2 } )$ & 0 & 0 & 0 & 70 & 595 & 2835
& 9975 \\ 
$ ( 2 1 ^{ 2 } )$ & 0 & 0 & 0 & 1260 & 16800 &
112770 & 525000 \\ 
$ ( 1 ^{ 4 } )$ & 0 & 0 & 0 & 2100 & 42000 &
389025 & 2368800 \\ 
\hline
$ ( 5 )$ & 0 & 0 & 0 & 1 & 8 & 36 & 120 \\ 
$ ( 4 1 )$ & 0 & 0 & 0 & 0 & 168 & 1596 &
8400 \\ 
$ ( 3 2 )$ & 0 & 0 & 0 & 0 & 280 & 2646 &
13860 \\ 
$ ( 3 1 ^{ 2 } )$ & 0 & 0 & 0 & 0 & 3360 & 47250
& 337680 \\ 
$ ( 2 ^{ 2 } 1 )$ & 0 & 0 & 0 & 0 & 5040 & 68040
& 474600 \\ 
$ ( 2 1 ^{ 3 } )$ & 0 & 0 & 0 & 0 & 47040 &
892080 & 8240400 \\ 
$ ( 1 ^{ 5 } )$ & 0 & 0 & 0 & 0 & 70560 &
1829520 & 22085280 \\ 
\hline
$ ( 6 )$ & 0 & 0 & 0 & 0 & 1 & 9 & 45 \\ 
$ ( 5 1 )$ & 0 & 0 & 0 & 0 & 0 & 252 & 2640
\\ 
$ ( 4 2 )$ & 0 & 0 & 0 & 0 & 0 & 504 & 5250
\\ 
$ ( 4 1 ^{ 2 } )$ & 0 & 0 & 0 & 0 & 0 & 7560
& 113400 \\ 
$ ( 3 ^{ 2 } )$ & 0 & 0 & 0 & 0 & 0 & 315 &
3276 \\ 
$ ( 3 2 1 )$ & 0 & 0 & 0 & 0 & 0 & 30240 &
428400 \\ 
$ ( 3 1 ^{ 3 } )$ & 0 & 0 & 0 & 0 & 0 & 181440
& 3477600 \\ 
$ ( 2 ^{ 3 } )$ & 0 & 0 & 0 & 0 & 0 & 7560 &
103600 \\ 
$ ( 2 ^{ 2 } 1 ^{ 2 } )$ & 0 & 0 & 0 & 0 & 0 &
430920 & 7862400 \\ 
$ ( 2 1 ^{ 4 } )$ & 0 & 0 & 0 & 0 & 0 & 2872800
& 68443200 \\ 
$ ( 1 ^{ 6 } )$ & 0 & 0 & 0 & 0 & 0 & 4011840
& 122673600 \\ 
\hline
$ ( 7 )$ & 0 & 0 & 0 & 0 & 0 & 1 & 10 \\ 
$ ( 6 1 )$ & 0 & 0 & 0 & 0 & 0 & 0 & 360 \\ 
$ ( 5 2 )$ & 0 & 0 & 0 & 0 & 0 & 0 & 840 \\ 
$ ( 5 1 ^{ 2 } )$ & 0 & 0 & 0 & 0 & 0 & 0 &
15120 \\ 
$ ( 4 3 )$ & 0 & 0 & 0 & 0 & 0 & 0 & 1260 \\

$ ( 4 2 1 )$ & 0 & 0 & 0 & 0 & 0 & 0 & 75600
\\ 
$ ( 4 1 ^{ 3 } )$ & 0 & 0 & 0 & 0 & 0 & 0 &
554400 \\ 
$ ( 3 ^{ 2 } 1 )$ & 0 & 0 & 0 & 0 & 0 & 0 &
50400 \\ 
$ ( 3 2 ^{ 2 } )$ & 0 & 0 & 0 & 0 & 0 & 0 &
75600 \\ 
$ ( 3 2 1 ^{ 2 } )$ & 0 & 0 & 0 & 0 & 0 & 0
& 3628800 \\ 
$ ( 3 1 ^{ 4 } )$ & 0 & 0 & 0 & 0 & 0 & 0 &
15120000 \\ 
$ ( 2 ^{ 3 } 1 )$ & 0 & 0 & 0 & 0 & 0 & 0 &
1890000 \\ 
$ ( 2 ^{ 2 } 1 ^{ 3 } )$ & 0 & 0 & 0 & 0 & 0 & 0
& 49140000 \\ 
$ ( 2 1 ^{ 5 } )$ & 0 & 0 & 0 & 0 & 0 & 0 &
264751200 \\ 
$ ( 1 ^{ 7 } )$ & 0 & 0 & 0 & 0 & 0 & 0 &
350481600
\\\hline
\end{tabular}
}
\label{tab1}
\end{table}
From Table \ref{tab1} we obtain the value of M\"obius function 
\begin{align*}
&\mu_{1+3,1}(\maxelm)=
3, \ 
\mu_{2+3,2}(\maxelm)=
-21, \ 
\mu_{3+3,3}(\maxelm)=
300, \ 
\mu_{4+3,4}(\maxelm)=
-7890, \\
&\mu_{5+3,5}(\maxelm)=
349650, \ 
\mu_{6+3,6}(\maxelm)=
-24188850,\ 
\mu_{7+3,7}(\maxelm)=2449878480,
\end{align*}
and characteristic polynomials
\allowdisplaybreaks
\begin{align*}
\chi_{1+3,1}=&
t^2 - 4t + 3 ,\\
\chi_{2+3,2}=&
t^3 - 10t^2 + 30t - 21, \\
\chi_{3+3,3}=&
t^4 - 20t^3 + 145t^2 - 426t + 300 ,\\
\chi_{4+3,4}=&
t^5 - 35t^4 + 490t^3 - 3381t^2 + 10815t - 7890 ,\\
\chi_{5+3,5}=&
t^6 - 56t^5 + 1330t^4 - 17136t^3 + 124971t^2 - 458760t + 349650 ,\\
\chi_{6+3,6}=&
t^7 - 84t^6 + 3108t^5 - 65926t^4 + 868014t^3 - 7039908t^2 \\
& +30423645t - 24188850, \\
\chi_{7+3,7}=&
t^8 - 120t^7 + 6510t^6 - 209692t^5 + 4414095t^4 - 62509140t^3
\\&+578334366t^2 - 2969914500t + 2449878.
\end{align*}

We can also compute the number of elements in $\VVV_2(l)$:
\begin{align*}
&\numof{\VVV_2(1)}= 1,\ 
\numof{\VVV_2(2)}= 2, \ 
\numof{\VVV_2(3)}=9, \
\numof{\VVV_2(4)}=96, \\
&\numof{\VVV_2(5)}=2419,\
\numof{\VVV_2(6)}= 133787, \ 
\numof{\VVV_2(7)}=14377347.
\end{align*}

\end{document}